\documentclass[10pt]{extarticle}

\usepackage[english]{babel}
\usepackage{graphicx}
\usepackage{framed}
\usepackage[normalem]{ulem}
\usepackage{indentfirst,ragged2e}
\usepackage{amsmath,amsthm,amssymb,amsfonts,wasysym,verbatim,bbm}
\usepackage{physics}
\usepackage[T1]{fontenc}
\usepackage{lmodern,mathrsfs}
\usepackage{wrapfig}
\usepackage[inline,shortlabels]{enumitem}
\setlist{topsep=2pt,itemsep=2pt,parsep=0pt,partopsep=0pt}
\usepackage[dvipsnames]{xcolor}
\usepackage[utf8]{inputenc}
\usepackage[a4paper,top=0.5in,bottom=0.2in,left=0.5in,right=0.5in,footskip=0.3in,includefoot]{geometry}
\usepackage{multicol}
\usepackage[most]{tcolorbox}
\usepackage{tikz,tikz-3dplot,tikz-cd,tkz-tab,tkz-euclide,tikzsymbols,pgf,pgfplots}
\pgfplotsset{compat=1.3}
\usepgfplotslibrary{fillbetween}
\usepackage{subfiles}
\graphicspath{{images/}{../images/}} 
\usepackage[backend=bibtex,style=numeric]{biblatex}
\usepackage{hyperref}
\usepackage[nameinlink]{cleveref}

\newcommand{\hide}[1]{} 

\definecolor{contcolor}{HTML}{514CE0}
\definecolor{convcolor}{HTML}{4CC6D4}

\newcommand{\ep}{\varepsilon}

\newcommand{\lam}{\lambda}

\newcommand{\N}{\mathbb{N}}

\newcommand{\R}{\mathbb{R}}

\newcommand{\As}{\mathcal{A}}
\newcommand{\Bs}{\mathcal{B}}

\newcommand{\Ps}{\mathcal{P}}

\newcommand{\Mf}{\mathfrak{M}}

\newcommand{\summ}[1][1]{\sum_{m=#1}^\infty}
\newcommand{\sumn}[1][1]{\sum_{n=#1}^\infty}

\newcommand{\emp}{\varnothing}

\newcommand{\sub}{\subseteq}

\newcommand{\kupp}{\bigsqcup}

\newcommand{\cappk}[1][1]{\bigcap_{k=#1}^\infty}
\newcommand{\cuppk}[1][1]{\bigcup_{k=#1}^\infty}

\newcommand{\kuppk}[1][1]{\bigsqcup_{k=#1}^\infty}
\newcommand{\kuppm}[1][1]{\bigsqcup_{m=#1}^\infty}
\newcommand{\kuppn}[1][1]{\bigsqcup_{n=#1}^\infty}

\newcommand{\pe}{\preceq}

\renewcommand{\iff}{\Leftrightarrow}

\newtheoremstyle{mystyle}{}{}{}{}{\sffamily\bfseries}{.}{ }{}
\makeatletter
\renewenvironment{proof}[1][\proofname] {\par\pushQED{\qed}{\normalfont\sffamily\bfseries\topsep6\p@\@plus6\p@\relax #1\@addpunct{.} }}{\popQED\endtrivlist\@endpefalse}
\makeatother

\theoremstyle{mystyle}{\newtheorem{definition}{Definition}}
\theoremstyle{mystyle}{}
\theoremstyle{mystyle}{\newtheorem{theorem}[definition]{Theorem}}
\theoremstyle{mystyle}{}
\theoremstyle{mystyle}{}
\theoremstyle{mystyle}{\newtheorem*{remark}{Remark}}
\theoremstyle{mystyle}{}
\theoremstyle{mystyle}{\newtheorem*{example}{Example}}
\theoremstyle{mystyle}{}
\theoremstyle{definition}{}

\newenvironment{talign*}{\csname align*\endcsname}{\endalign}

\DeclareMathAlphabet\mathbfcal{OMS}{cmsy}{b}{n}
\newcommand{\ind}{\hspace{0.2in}} 

\makeatletter
\g@addto@macro\normalsize{
\setlength\abovedisplayskip{3pt}
\setlength\belowdisplayskip{3pt}
\setlength\abovedisplayshortskip{0pt}
\setlength\belowdisplayshortskip{0pt}}
\makeatother

\makeatletter
\renewenvironment{abstract}{
\if@twocolumn
\section*{\abstractname}
\else
\begin{center}
{\sffamily\bfseries\abstractname\vspace{\z@}}
\end{center}
\quotation
\fi}
{\if@twocolumn\else\endquotation\fi}
\makeatother

\title{\Large\sffamily\bfseries Measures Form a Complete Lattice}
\author{\Large\sffamily Senan Sekhon}
\date{\large\sffamily April 14, 2021}

\begin{document}

\maketitle

\begin{abstract}
    We show that the set of all measures on any measurable space is a complete lattice, i.e. every collection of measures has both a greatest lower bound and a least upper bound.
\end{abstract}

\ind In a first course on measure theory, one learns that the sum of any collection of measures, as well as any non-negative constant multiple of a measure, is also a measure on the same measurable space. Additionally, the limit of an increasing sequence of measures is also a measure on the same measurable space.\\

\ind On the other hand, the minimum and maximum of two measures are not measures (in general). This is shown by the following example:

\begin{example}
Suppose $X=\{a,b\}$, $\As=\Ps(X)$ and define $\mu,\nu:\As\to[0,\infty]$ as follows:
\begin{align*}
    \mu(\emp)=0 && \mu(\{a\})=1 && \mu(\{b\})=0 && \mu(\{a,b\})=1 \\
    \nu(\emp)=0 && \nu(\{a\})=0 && \nu(\{b\})=1 && \nu(\{a,b\})=1
\end{align*}
In other words, $\mu$ and $\nu$ are the Dirac measures at $a$ and $b$ respectively. Then we have:
\begin{align*}
    \min\{\mu,\nu\}(\{a\})+\min\{\mu,\nu\}(\{b\})=0+0=0&\ne1=\min\{\mu,\nu\}(\{a,b\}) \\
    \max\{\mu,\nu\}(\{a\})+\max\{\mu,\nu\}(\{b\})=1+1=2&\ne1=\max\{\mu,\nu\}(\{a,b\})
\end{align*}
Thus $\min\{\mu,\nu\}$ and $\max\{\mu,\nu\}$ are not measures on $(X,\As)$.
\end{example}

\ind Thus, if we want to establish a ``greatest lower bound'' or a ``least upper bound'' of two measures, we need to be more careful. We first define precisely the partial order in use:

\begin{definition}
Suppose $(X,\As)$ is a measurable space. Then $\Mf(X,\As)$ denotes the poset of all measures on $(X,\As)$, with the following partial order:
\begin{equation*}
    \mu\pe\nu \iff \mu(A)\le\nu(A) \text{ for all } A\in\As
\end{equation*}
\end{definition}
\ind It follows directly from the definition that $\pe$ is a partial order on $\Mf(X,\As)$. From now on, we will denote $\Mf(X,\As)$ by $\Mf$, to indicate that the measurable space is arbitrary.\\

\ind Clearly $\Mf$ has a least element, the \emph{zero measure} ($\mu(A)=0$ for all $A\in\As$). It also has a greatest element, the \emph{infinity measure} ($\mu(\emp)=0$ and $\mu(A)=\infty$ for all $A\in\As$, $A\ne\emp$).\\

\ind Going back to the question of the greatest lower bound and least upper bound, we will present two ways to motivate the answer:\\

{\sffamily\bfseries Method 1.} Suppose $\lam\in\Mf$ such that $\lam\pe\mu$ and $\lam\pe\nu$. Then for all $A\in\As$, we must have $\lam(A)\le\mu(A)$ and $\lam(A)\le\nu(A)$. This alone only implies that $\lam(A)\le\min\{\mu(A),\nu(A)\}$, which we have already seen is not good enough. To optimize this, we can introduce another set $B\in\As$ and do the following:
\begin{equation*}
    \lam(A)=\lam(A\cap B)+\lam(A\cap B^c)\le\mu(A\cap B)+\nu(A\cap B^c)
\end{equation*}
In other words, we must have $\lam(A)\le\mu(A\cap B)+\nu(A\cap B^c)$ for \emph{every} $B\in\As$. This motivates the following definition:
\begin{equation*}
    \lam(A)=\inf_{B\in\As} \{\mu(A\cap B)+\nu(A\cap B^c)\}
\end{equation*}
In the example of the Dirac measures, this definition works, as it produces the zero measure.\\

{\sffamily\bfseries Method 2.} If $f:S\to\R$ is a real-valued function (on any set $S$), we can define its \emph{positive part} $f^+=\max\{f,0\}$ and \emph{negative part} $f^-=-\min\{f,0\}$. These form a decomposition of $f$, since $f=f^+-f^-$. Also $f^+$ and $f^-$ are `complementary', in the sense that at every point $x\in S$, at least one of $f^+(x)$ and $f^-(x)=0$.\\

\ind In the case of signed measures, there is an analogous notion, known as \emph{Jordan decomposition}: If $\lam$ is a signed measure, then we can write $\lam=\lam^+-\lam^-$, where $\lam^+$ and $\lam^-$ are measures and $\lam^+\perp\lam^-$. See \cite[Theorem 9.30, Page 269]{axler} or \cite[Theorem 3.4, Page 87]{folland} for a proof.\\

\ind If $f$ and $g$ are real-valued functions, we can express $\min\{f,g\}$ and $\max\{f,g\}$ using positive and negative parts:
\begin{align*}
    \min\{f,g\}=f+\min\{0,g-f\}=f-(g-f)^- && \max\{f,g\}=f+\max\{0,g-f\}=f+(g-f)^+
\end{align*}
This motivates the following definitions:
\begin{align*}
    \mu\wedge\nu=\mu-(\nu-\mu)^- && \mu\vee\nu=\mu+(\nu-\mu)^+
\end{align*}
However, this is not always valid. In particular, $\nu-\mu$ may not be well-defined (due to $\infty-\infty$). When it is valid, though, it is equivalent to the definition from Method 1, which can be seen using the following formulas:
\begin{align*}
    \lam^+(A)=\sup\{\lam(E)\mid E\in\As,E\sub A\} && \lam^-(A)=-\inf\{\lam(E)\mid E\in\As,E\sub A\}
\end{align*}

\vspace{10pt}

\ind We are now ready to answer the question of the greatest lower bound and least upper bound of two measures. From now on, we will use $\sqcup$ to denote \emph{disjoint} unions, and $\cup$ to denote arbitrary unions.

\begin{theorem}\label{glblub}
Suppose $\mu,\nu\in\Mf$. Then their greatest lower bound $\mu\wedge\nu$ and least upper bound $\mu\vee\nu$ are given by:
\begin{align*}
    (\mu\wedge\nu)(A)=\inf_{B\in\As} \left\{\mu(A\cap B)+\nu(A\cap B^c)\right\} &&
    (\mu\vee\nu)(A)=\sup_{B\in\As} \left\{\mu(A\cap B)+\nu(A\cap B^c)\right\}
\end{align*}
\end{theorem}
\begin{proof}
We will prove the result for $\mu\wedge\nu$, the result for $\mu\vee\nu$ can be proved similarly. Define $\lam=\mu\wedge\nu$. We first show that $\lam\in\Mf$, i.e. $\lam$ is a measure on $(X,\As)$.\\

For all $B\in\As$, we have $\mu(\emp\cap B)+\nu(\emp\cap B^c)=0+0=0$, so $\lam(\emp)=0$.\\

Suppose $\{A_n\}_{n=1}^\infty$ is a countable collection of disjoint sets in $\As$, and define $A=\kuppn A_n$. Then for all $B\in\As$, we have:
\begin{align*}
    A\cap B=\left(\kuppn A_n\right)\cap B=\kuppn(A_n\cap B) && A\cap B^c=\left(\kuppn A_n\right)\cap B^c=\kuppn(A_n\cap B^c)
\end{align*}
This yields:
\begin{align*}
    \mu(A\cap B)+\nu(A\cap B^c)
    &=\mu\left(\kuppn(A_n\cap B)\right)+\nu\left(\kuppn(A_n\cap B^c)\right)
    =\sumn\mu(A_n\cap B)+\sumn\nu(A_n\cap B^c)\\
    &=\sumn(\mu(A_n\cap B)+\nu(A_n\cap B^c))
    \ge\sumn\lam(A_n)
\end{align*}
Taking the infimum over all $B\in\As$ yields $\lam(A)\ge\sumn\lam(A_n)$.\\

If $\sumn\lam(A_n)=\infty$, then we have $\lam(A)\ge\sumn\lam(A_n)$ as both sides are $\infty$. Suppose $\sumn\lam(A_n)<\infty$. Then for any $n\in\N$ and any $\ep>0$, there is a set $B_n\in\As$ such that:
\begin{equation}\label{ep1}
    \mu(A_n\cap B_n)+\nu(A_n\cap B_n^c)<\lam(A_n)+\frac{\ep}{2^n}
\end{equation}
Define $C_n=A_n\cap B_n$ and $C=\kuppn C_n$ (this is a disjoint union as $C_n\sub A_n$ and all $A_n$ are disjoint). Since $A_n,B_n\in\As$ for all $n\in\N$, so are all $C_n$ and $C$. Note that:
\begin{align}
    A_n\cap C&=A_n\cap\left(\kuppk C_k\right)=\cuppk(A_n\cap C_k)=A_n\cap C_n \tag{since $A_n\cap C_k=\emp$ for $k\ne n$}\\
    &=A_n\cap(A_n\cap B_n)=A_n\cap B_n \label{abc1}\\
    A_n\cap C^c&=A_n\cap\left(\kuppk C_k\right)^c=A_n\cap\left(\cappk C_k^c\right)=\cappk(A_n\cap C_k^c)=A_n\cap C_n^c \tag{since $A_n\cap C_k^c=A_n$ for $k\ne n$}\\
    &=A_n\cap(A_n\cap B_n)^c=A_n\cap(A_n^c\cup B_n^c)=A_n\cap B_n^c \label{abc2}
\end{align}
This yields:
\begin{align*}
    \mu(A\cap C)+\nu(A\cap C^c)
    &=\mu\left(\kuppn(A_n\cap C)\right)+\nu\left(\kuppn(A_n\cap C^c)\right)
    =\sumn\mu(A_n\cap C)+\sumn\nu(A_n\cap C^c) \\
    &=\sumn(\mu(A_n\cap C)+\nu(A_n\cap C^c)) \\
    &=\sumn(\mu(A_n\cap B_n)+\nu(A_n\cap B_n^c)) \tag{by \eqref{abc1} and \eqref{abc2}}\\
    &<\sumn\left(\lam(A_n)+\frac{\ep}{2^n}\right) \tag{by \eqref{ep1}}\\
    &=\sumn\lam(A_n)+\sumn\frac{\ep}{2^n}
    =\sumn\lam(A_n)+\ep
\end{align*}
Thus $\lam(A)<\sumn\lam(A_n)+\ep$. Since this holds for all $\ep>0$, we have $\lam(A)\le\sumn\lam(A_n)$, and so:
\begin{equation*}
    \lam(A)=\sumn\lam(A_n)
\end{equation*}
Thus $\lam$ is a measure on $(X,\As)$.\\

For all $A\in\As$, we have $\lam(A)\le\mu(A)$ and $\lam(A)\le\nu(A)$ (these follow by setting $B=X$ and $B=\emp$ respectively in the definition). Thus $\lam\pe\nu$ and $\lam\pe\nu$. Suppose $\rho\in\Mf$ such that $\rho\pe\mu$ and $\rho\pe\nu$. Then for all $A,B\in\As$, we have:
\begin{equation*}
    \rho(A)=\rho(A\cap B)+\rho(A\cap B^c)\le\mu(A\cap B)+\nu(A\cap B^c)
\end{equation*}
Taking the infimum over all $B\in\As$ yields $\rho(A)\le\lam(A)$. Thus $\lam$ is the greatest lower bound of $\mu$ and $\nu$.
\end{proof}

\vspace{10pt}

\ind \Cref{glblub} appears as an exercise in \cite[Exercise 9B.8, Page 278]{axler}, \cite[Exercise 3.1.7a, Page 88]{folland} and \cite[Exercise 6.76, Page 364]{mcdonald}. The formulas also appear in \cite[Ch. III, §1, No. 5, Th. 3,Page III.12]{bourbaki} and \cite[Page 111]{dales} (again, without proof), though only for a special class of Borel measures on a locally compact topological space.\\

\ind This proves that $\Mf$ is a lattice, i.e. a poset in which every pair of elements has a greatest lower bound and a least upper bound. We will show in \Cref{glblubi} that it is a \emph{complete} lattice, i.e. \emph{every} (possibly uncountable) subset of $\Mf$ has a greatest lower bound and a least upper bound.\\

\ind To motivate the answer, let's go back to the formula $(\mu\wedge\nu)(A)=\inf_{B\in\As} \left\{\mu(A\cap B)+\nu(A\cap B^c)\right\}$. Clearly we can freely switch $B$ and $B^c$, since $\As$ is a $\sigma$-algebra and is thus closed under complements. The important bit here is that $B$ and $B^c$ \emph{partition} $X$. In other words, we split $X$ into two pieces, give one to $\mu$ and the other to $\nu$. Of course, this readily generalizes to arbitrary collections of sets. With this, you might think of the following definition of the greatest lower bound of a collection $\{\mu_\alpha\}_{\alpha\in I}\sub\Mf$:
\begin{equation}\label{wrongd}
    \lam(A)=\inf\left\{\sum_{\alpha\in I} \mu_\alpha(A\cap B_\alpha)\right\}
\end{equation}
Where the infimum is over all partitions $\{B_\alpha\}_{\alpha\in I}$ of $X$ into measurable sets $B_\alpha\in\As$. While this does yield a lower bound of $\{\mu_\alpha\}_{\alpha\in I}$, it may not be the \emph{greatest} lower bound, as the following example shows.

\begin{example}
Suppose $X=[0,1]$, $\As=\Bs([0,1])$ (the Borel $\sigma$-algebra on $[0,1]$) and $I=[4,5]$. For each $\alpha\in[4,5]$, define $\mu_\alpha=\alpha\mu_L$ ($\alpha$ times the Lebesgue measure). If we partition $[0,1]$ into individual points, say $B_\alpha=\{\alpha-4\}$, then we have:
\begin{equation*}
    \sum_{\alpha\in[4,5]} \mu_\alpha(\{\alpha-4\})=\sum_{\alpha\in[4,5]} 0=0
\end{equation*}
Thus $\lam([0,1])=0$, i.e. $\lam$ is the zero measure. However, for all $\alpha\in[4,5]$, we have $\mu_\alpha\ge4\mu_L$, so $4\mu_L$ is a greater lower bound.
\end{example}

\begin{remark}
The problem with \eqref{wrongd} is not that the sum may be uncountable: \emph{any} sum of non-negative extended real numbers is well-defined as an extended real number. The problem is that measures need not be additive over uncountable collections of sets, so the formula $\mu(\kupp_{n\in\N} A_n)=\sum_{n\in\N} \mu(A_n)$ does not apply if we replace $\N$ with an uncountable set.
\end{remark}

\ind The above example shows that \eqref{wrongd} is not the correct formula. Before stating what the correct formula is, we make the following definition:

\begin{definition}
Suppose $(X,\As)$ is a measurable space and $\{\mu_\alpha\}_{\alpha\in I}$ is a collection of measures on $(X,\As)$. Then $\Sigma$ denotes the set of all pairs $(\{B_n\}_{n=1}^\infty,\{\mu_{\alpha_n}\}_{n=1}^\infty)$, where $\{B_n\}_{n=1}^\infty$ is a countable collection of disjoint sets in $\As$ such that $\kuppn B_n=X$ (also known as a \emph{countable measurable partition} of $X$), and for each $n\in\N$, $\mu_{\alpha_n}$ is some measure in the collection $\{\mu_\alpha\}_{\alpha\in I}$.
\end{definition}
In other words, we choose a countable\footnote{We will assume without loss of generality that it is countably infinite, since if it is finite, we can extend it by making the rest of the sets empty.} collection of measurable sets that partition $X$, and for each of these sets, we attach one of the measures\footnote{We will assume without loss of generality that all of these measures are distinct, since if any of them coincide, we can merge the corresponding sets from the partition.} from the collection $\{\mu_\alpha\}_{\alpha\in I}$.\\

In the proof of \Cref{glblubi}, we will use the fact that if $a_{m,n}\ge0$ for all $m,n\in\N$, then:
\begin{equation}\label{tcm}
    \summ\sumn a_{m,n}=\sumn\summ a_{m,n}
\end{equation}
This is simply Tonelli's theorem applied to the counting measure on $(\N,\Ps(\N))$. See \cite[Theorem 5.31, Page 131]{axler} or \cite[Example 4.23b, Page 249]{mcdonald} for a proof.

\begin{theorem}\label{glblubi}
Suppose $\{\mu_\alpha\}_{\alpha\in I}$ is a collection of measures on $(X,\As)$. Then their greatest lower bound $\bigwedge_{\alpha\in I} \mu_\alpha$ and least upper bound $\bigvee_{\alpha\in I} \mu_\alpha$ are given by:
\begin{align*}
    \left(\bigwedge_{\alpha\in I} \mu_\alpha\right)(A)=\inf\left\{\sumn\mu_{\alpha_n}(A\cap B_n)\right\} && \left(\bigvee_{\alpha\in I} \mu_\alpha\right)(A)=\sup\left\{\sumn\mu_{\alpha_n}(A\cap B_n)\right\}
\end{align*}
Where the infimum and supremum are over all $(\{B_n\}_{n=1}^\infty,\{\mu_{\alpha_n}\}_{n=1}^\infty)\in\Sigma$.
\end{theorem}
\begin{proof}
We will prove the result for $\bigwedge_{\alpha\in I} \mu_\alpha$, the result for $\bigvee_{\alpha\in I} \mu_\alpha$ can be proved similarly. Define $\lam=\bigwedge_{\alpha\in I} \mu_\alpha$. We will use the shorthand $\{B_n,\mu_{\alpha_n}\}$ for a pair $(\{B_n\}_{n=1}^\infty,\{\mu_{\alpha_n}\}_{n=1}^\infty)\in\Sigma$.\\

We first show that $\lam\in\Mf$, i.e. $\lam$ is a measure on $(X,\As)$.\\

For all measurable partitions $\{B_n\}_{n=1}^\infty$ of $X$, we have $\sumn\mu_{\alpha_n}(\emp\cap B_n)=\sumn\mu_{\alpha_n}(\emp)=\sumn 0=0$, so $\lam(\emp)=0$.\\

Suppose $\{A_m\}_{m=1}^\infty$ is a countable collection of disjoint sets in $\As$, and define $A=\kuppm A_m$. Then for any $\{B_n,\mu_{\alpha_n}\}$, we have:
\begin{equation*}
    A\cap B_n=\left(\kuppm A_m\right)\cap B_n=\kuppm(A_m\cap B_n)
\end{equation*}
This yields:
\begin{align*}
    \sumn\mu_{\alpha_n}(A\cap B_n)
    &=\sumn\mu_{\alpha_n}\left(\kuppm(A_m\cap B_n)\right)
    =\sumn\summ\mu_{\alpha_n}(A_m\cap B_n) \\
    &=\summ\sumn\mu_{\alpha_n}(A_m\cap B_n)
    \ge\summ\lam(A_m) \tag{by \eqref{tcm}}
\end{align*}
Taking the infimum over all $\{B_n,\mu_{\alpha_n}\}$ yields $\lam(A)\ge\summ\lam(A_m)$.\\

If $\summ\lam(A_m)=\infty$, then we have $\lam(A)=\summ\lam(A_m)$ as both sides are $\infty$. Suppose $\summ\lam(A_m)<\infty$. Then for any $m\in\N$ and any $\ep>0$, there exists $\{B_n,\mu_{\alpha_n}\}$ such that:
\begin{equation}\label{ep2}
    \sumn\mu_{\alpha_n}(A_m\cap B_{m,n})<\lam(A_m)+\frac{\ep}{2^m}
\end{equation}
Define $C_{m,n}=A_m\cap B_{m,n}$ and $C_n=\kuppm C_{m,n}$ (this is a disjoint union as $C_{m,n}\sub A_m$ and all $A_m$ are disjoint). Since $A_m,B_{m,n}\in\As$ for all $m,n\in\N$, so are all $C_{m,n}$ and all $C_n$. We also have:
\begin{align}
    A_m\cap C_n&=A_m\cap\left(\kuppk C_{k,n}\right)=\kuppk(A_m\cap C_{k,n})=A_m\cap C_{m,n} \tag{since $A_m\cap C_{k,n}=\emp$ for $k\ne m$}\\
    &=A_m\cap(A_m\cap B_{m,n})=A_m\cap B_{m,n} \label{abcn}
\end{align}
This yields:
\begin{align*}
    \sumn\mu_{\alpha_n}(A\cap C_n)
    &=\sumn\mu_{\alpha_n}\left(\kuppm(A_m\cap C_n)\right)
    =\sumn\summ\mu_{\alpha_n}(A_m\cap C_n) \\
    &=\summ\sumn\mu_{\alpha_n}(A_m\cap C_n)  \tag{by \eqref{tcm}}\\
    &=\summ\sumn\mu_{\alpha_n}(A_m\cap B_{m,n}) \tag{by \eqref{abcn}}\\
    &<\summ\left(\lam(A_m)+\frac{\ep}{2^m}\right) \tag{by \eqref{ep2}}\\
    &=\summ\lam(A_m)+\summ\frac{\ep}{2^m}
    =\summ\lam(A_m)+\ep
\end{align*}
Again, we can swap the sums as all summands are non-negative. While the sets $C_n$ are disjoint, they do not partition $X$ (since $\kuppn C_n=A$, not $X$). We can fix this by redefining one of the sets $C_k$ to be $C_k\sqcup A^c$. This does not affect the steps above, as we only considered the measures of intersections with $A_m$ or $A$ (and so they are subsets of $A$). With this, $\{C_n\}_{n=1}^\infty$ becomes a countable measurable partition of $X$ such that:
\begin{equation*}
    \sumn\mu_{\alpha_n}(A\cap C_n)<\summ\lam(A_m)+\ep
\end{equation*}
Thus $\lam(A)<\summ\lam(A_m)+\ep$. Since this holds for all $\ep>0$, we have $\lam(A)\le\summ\lam(A_m)$, and so:
\begin{equation*}
    \lam(A)=\summ\lam(A_m)
\end{equation*}
Thus $\lam$ is a measure on $(X,\As)$.\\

For all $\alpha\in I$ and all $A\in\As$, we have $\lam(A)\le\mu_\alpha(A)$ (this follows by setting $\alpha_1=\alpha$, $B_{\alpha_1}=X$ and $B_{\alpha_n}=\emp$ for all other $n>1$ in the definition of $\lam$). Suppose $\rho\in\Mf$ such that $\rho\pe\mu_\alpha$ for all $\alpha\in I$. Then for all $A\in\As$ and all $\{B_n,\mu_{\alpha_n}\}$, we have:
\begin{equation*}
    \rho(A)=\sumn\rho(A\cap B_n)\le\sumn\mu_{\alpha_n}(A\cap B_n)
\end{equation*}
Taking the infimum over all $\{B_n,\mu_{\alpha_n}\}$ yields $\rho(A)\le\lam(A)$. Thus $\lam$ is the greatest lower bound of $\{\mu_\alpha\}_{\alpha\in I}$.
\end{proof}

\printbibliography

\end{document}